\documentclass[12pt, a4paper]{amsart}
\usepackage{amscd, amssymb}
\usepackage{amsthm}
\def\C{\mathbb{C}}
\def\newspan{\operatorname{span}}
\newtheorem{thm}{Theorem}

\begin{document}
\title[Every AF-algebra is Morita equivalent to a graph algebra.]{Every
AF-algebra is Morita equivalent\\
to a graph algebra.}

\author{Jason Tyler}

\address{School of Mathematical and Physical Sciences\\
University of Newcastle\\
NSW 2308\\
Australia}
\email{jason.tyler@studentmail.newcastle.edu.au}

\date{3 June 2003}

\thanks{My thanks go to my indefatigable supervisor Iain Raeburn for
his guidance throughout this work.}

\begin{abstract}
We show how to modify any Bratteli diagram $E$ for an AF-algebra $A$ to
obtain a Bratteli diagram $KE$ for $A$ whose graph algebra $C^*(KE)$
contains both $A$ and $C^*(E)$ as full corners.
\end{abstract}

\maketitle

An elegant theorem of Drinen says that every AF-algebra $A$ is isomorphic
to a corner in a graph algebra~\cite[Theorem 1]{dri}, and hence is Morita
equivalent to the graph algebra.
The graph in question is a Bratteli diagram for $A$,
but it needs to be a carefully chosen one; two constructions of such a
diagram were
described in~\cite{dri}, one attributed to Kumjian. Here we show that
applying Kumjian's construction to an arbitrary Bratteli diagram $E$
for $A$ gives a graph $KE$ whose $C^*$-algebra contains both $A$ and
$C^*(E)$ as full corners, so that $A$ is Morita equivalent to
the $C^*$-algebra $C^*(E)$ of the original Bratteli diagram $E$.

A \emph{directed graph} $E$ consists of
countable sets $E^0$ of vertices and $E^1$ of edges, along with functions
$r,s:E^1\to E^0$ which map edges to their range and source vertices.
The graph is \emph{row-finite} if each vertex emits at most finitely
many edges.
Given a row-finite graph $E$, a \emph{Cuntz-Krieger $E$-family} in a
$C^*$-algebra consists of a set of mutually orthogonal projections
$\{p_v:v\in E^0\}$ and a set of partial isometries $\{s_e:e\in E^1\}$
satisfying the \emph{Cuntz-Krieger relations}:
\[
s_e^*s_e = p_{r(e)} \text{ for }e\in E^1\text{ and }
p_v=\sum_{e\in s^{-1}(v)}s_es_e^*\text{ whenever }
s^{-1}(v)\not=\emptyset.
\]
The graph algebra $C^*(E)$ is the universal $C^*$-algebra generated by a
Cuntz-Krieger $E$-family $\{s_e,p_v\}$~\cite[Theorem 1.2]{kpr}. We denote
by $E^*$ the set of all finite
paths in $E$; that is, sequences of
edges $\mu_1\mu_2\ldots\mu_n$ such that $r(\mu_i)=s(\mu_{i+1})$ for 
$1\le i<n$.  We include the vertices as paths of length zero.
Given $\mu=\mu_1\mu_2\ldots\mu_n\in E^*$, define $s_\mu:=
s_{\mu_1}s_{\mu_2}\ldots s_{\mu_n}$. 
It follows from~\cite[Lemma 1.1]{kpr}
that \[
C^*(E) =\overline{\newspan}\{s_\mu s_\nu^*: \mu,\nu\in E^*,
r(\mu)=r(\nu)\}.
\]

A \emph{Bratteli diagram} is
a  directed graph $E$ such that:

\begin{itemize}

\item $E^0$ is the disjoint union of finite sets $\{V_n\}$,

\item every edge with source in $V_n$ has range in $V_{n+1}$, and

\item each $v\in E^0$ is labelled with a positive
integer $d_v$ satisfying $d_v\ge \sum_{e\in r^{-1}(v)} d_{s(e)}$.

\end{itemize}
We say that $E$ is a Bratteli diagram for a sequence
of $C^*$-algebras $A_1\subset A_2\subset \ldots$ if each $A_n$ is
isomorphic to 
$\bigoplus_{v\in V_n} M_{d_v}(\C)$ and
the embedding of each $M_{d_v}(\C)\subset A_n$ in each $M_{d_w}(\C)
\subset A_{n+1}$  scales the trace by $\#(s^{-1}(v)\cap r^{-1}(w))$.
We say that $E$ is a Bratteli diagram for an AF-algebra $A$ if there
exists a sequence of $C^*$-subalgebras $\{A_n\}$ of $A$ such that
$A=\overline{\bigcup A_n}$ and $E$ is a Bratteli diagram for $\{A_n\}$.

\begin{thm}
Let $E$ be a Bratteli diagram for an AF-algebra $A$.
Then there exists a Bratteli diagram $KE$ for $A$ such that
$C^*(KE)$ contains $A$ and $C^*(E)$ as complementary full corners.
\end{thm}

The projection $p$ defining the corner is the sum $p=\sum_{v\in S} p_v $
where $S\subset KE^0$; this sum converges strictly to a projection
in $M(C^*(KE))$ by~\cite[Lemma 1.1]{bprs}. Crucial for us is the
observation that for $\mu,\nu\in KE^*$, \[
ps_\mu s_\nu^*=\begin{cases}
s_\mu s_\nu^*&\text{if $s(\mu)\in S$}\\
0&\text{otherwise}
\end{cases}
\] so that \[
pC^*(KE)p=\overline{\newspan}\{s_\mu s_\nu^*: s(\mu),s(\nu)\in S,
r(\mu)=r(\nu)\}.
\]

\begin{proof}[Proof of the theorem.]
For $n>0$, denote by $V_n$ the set of vertices on the $n$th level of
$E$, and let $V_0=\emptyset$.
For each $v\in E^0$, let $d_v$ be the rank
of the matrix algebra corresponding to $v$. 
For every vertex
$v\in E^0$, calculate $\sigma_v:=d_v - 
\sum_{e\in r^{-1}(v)}d_{s(e)}$.
We define $KE^0=\bigcup_{n=0}^\infty KV_n$, where \[
KV_n := \begin{cases}
V_n & \text{if $\sigma_v=0$ for all $v\in V_{n+1}$}\\
V_n \cup \{w_n\} & \text{if $\sigma_v>0$ for some $v\in V_{n+1}$,}
\end{cases} \]
and define $KE^1$ to be $E^1$ together with, for every $w_n$ and 
$v\in V_{n+1}$,
$\sigma_v$ edges from $w_n$ to $v$.
Denote by $S$ the set $KE^0\backslash E^0=\bigcup \{w_n\}$,
and set $d_w=1$ for all $w\in S$.
Constructing $KE$ in this fashion ensures
that for all $v\in {KE}^0$, the number of paths beginning in $S$
and ending at $v$ is $d_v$. Note that if $A$ is unital,
so $\sigma_v=0$ for all $v\in E^0\backslash V_1$, then we add only one
vertex to $E$ and $KE^0 = E^0 \bigcup \{w_0\}$.

Since $E$ is a Bratteli diagram for $A$, there is an increasing sequence of 
$C^*$-subalgebras $F_n$ of $A$ such that $A=\overline{\bigcup F_n}$ and
$E$ is a Bratteli diagram for the sequence $\{F_n\}$.
For those $n$ where $KV_n \not= V_n$,
we define a subalgebra $F'_n$ of $A$ by $F'_0:= \C1$ and
\[
F'_n:=F_n\oplus\C(1_{F_{n+1}}-1_{F_n})\cong\bigoplus_{v\in
V_n} M_{d_v}(\C)\oplus \C\text{\hskip1cm for $n>0$.} \]
For all other $n$, define $F'_n=F_n$. The graph $KE$ is then a
Bratteli diagram for the sequence $\{F'_n\}$.  Since $F_n\subseteq
F'_n\subseteq F_{n+1}$ for all $n$, we have $\overline{\bigcup F'_n} =
\overline {\bigcup F_n} = A$; thus $KE$ is a Bratteli diagram for $A$.

Let $\{s_e,p_v\}$ be the universal Cuntz-Krieger $KE$-family generating
$C^*(KE)$. Define a projection 
$p\in M(C^*(KE))$ by $p:=\sum_{v\in S}p_v$.
We aim to show that the corner $pC^*(KE)p$
is isomorphic to $A$. Since two algebras with the same
Bratteli diagram are isomorphic~\cite[Proposition III.2.7]{dav}, we can
achieve this by identifying a sequence of subalgebras of $pC^*(KE)p$
for which $E$ is a Bratteli diagram and whose union is dense in $pC^*(KE)p$.
For each $n>0$ define $D_n:= \newspan\{D^v:v\in V_n\}$,
where \[
D^v:= \newspan\{s_\mu s_\nu^*: \mu,\nu\in KE^*,s(\mu),s(\nu)\in S,
r(\mu)=r(\nu)=v\}
\] for each $v\in KE^0$.
Note that
\[
pC^*(KE)p =\overline{\newspan}\{s_\mu s_\nu^*:
\mu,\nu\in KE^*,s(\mu),
s(\nu)\in S,r(\mu)=r(\nu)\} = \overline{\bigcup D_n}.\]
Given $v\in E^0$ and paths $\mu,\nu,\alpha,\beta$ with source in $S$ and
range $v$, observe that none of $\mu,\nu,\alpha,\beta$ can extend any
other since $KE$ contains no loops; \cite[Lemma 1.1]{kpr} then gives \[
s_\mu s_\nu^*s_\alpha s_\beta^* = \begin{cases}
s_\mu s_\beta^*&\text{if $\nu=\alpha$}\\
0 &\text{otherwise.}
\end{cases}
\]
Also, $(s_\mu s_\nu^*)^* = s_\nu s_\mu^*$, so
\[\{s_\mu s_\nu^*: \mu,\nu\in KE^*,s(\mu),s(\nu)\in S, r(\mu)=r(\nu)=v\}\]
is a family of matrix units. Since there are $d_v$ paths $\mu$ with
$s(\mu)\in S$ and  $r(\mu)=v$, $D^v$ is isomorphic to $M_{d_v}(\C)$.
Further, note that for distinct $v,w\in V_n$, no path ending at $v$
may extend one ending at $w$, so $D^vD^w=0$ and $D_n = \bigoplus_{v\in
V_n}D^v\cong \bigoplus_{v\in V_n} M_{d_v}(\C)$. It remains only
to check that the embedding of each $D_n$ in $D_{n+1}$ matches that
described by $E$; specifically, for $v\in V_n$ and $w\in
V_{n+1}$ we need that $D^v$ is embedded in $D^w$ with multiplicity
$\#(s^{-1}(v)\cap r^{-1}(w))$.
This follows from the Cuntz-Krieger relations at
$v$: take paths $\mu,\nu$ with source in $S$ and range $v$, decompose
the matrix unit $s_\mu s_\nu^*\in D^v$ as \[
s_\mu s_\nu^*=s_\mu p_v s_\nu^*
=s_\mu \Bigl(\sum_{e\in s^{-1}(v)}s_e s_e^*\Bigr) s_\nu^*
=\sum_{e\in s^{-1}(v)}s_{\mu e}s_{\nu e}^* \]
and note that $s_{\mu e}s_{\nu e}^* $ is a matrix unit in $D^w$
precisely when $e\in r^{-1}(w)$.

Consider now the complementary corner \[
(1-p)C^*(KE)(1-p)=
\overline{\newspan}\{s_\mu s_\nu^*:\mu,\nu\in KE^*,
s(\mu),s(\nu)\in E^0,r(\mu)=r(\nu)\}.\]
Since $KE^1 \backslash E^1$ contains only edges from $S$ to $E^0$,
paths beginning in $E^0$ never leave $E^0$. Thus $(1-p)C^*(KE)(1-p)$ is 
generated by the Cuntz-Krieger $E$-family \[\{s_e,p_v:e\in E^1,v\in E^0\}.\]
Further, $E$ contains no loops, so the Cuntz-Krieger uniqueness 
theorem~\cite[Theorem 3.1]{bprs} implies that $(1-p)C^*(KE)(1-p)$ is 
isomorphic to $C^*(E)$.

Finally, we must show that $p$ and $1-p$ are full. 
Note that for every $v\in KE^0$ there is a
path beginning in $S$ and ending at $v$.
Suppose that $I$ is an ideal in $C^*(KE)$
containing $pC^*(KE)p$; then $I$ certainly contains the projections
$\{p_w:w\in S\}$. 
Given a vertex $v$ in $E^0$, choose a path $\alpha$ beginning at some $w\in S$
and ending at $v$.
Then $s_\alpha= p_{w}s_\alpha\in I$, so $p_v=s_\alpha^*s_\alpha\in
I$, every generator $\{s_e,p_v\}$ of $C^*(KE)$ is in $I$, and
$I=C^*(KE)$. Now suppose that $J$ is an ideal in $C^*(KE)$ containing
$(1-p)C^*(KE)(1-p)$, so for every $v\in E^0$ we have $p_v\in J$. Given a
vertex $v\in S$, note that every edge $e$ with $s(e)=v$ satisfies
$r(e)\in E^0$; so for all $e\in s^{-1}(v)$,
we know that $p_{r(e)}=s_e^*s_e\in J$, implying $s_e= s_es_e^*s_e\in
J$ and $s_es_e^*\in J$. Thus $p_v=\sum_{e\in s^{-1}(v)}s_es_e^*\in J$,
the universal $KE$-family $\{s_e,p_v\}$ is contained in $J$,
and $J=C^*(KE)$.
\end{proof}

\end{document}